\documentclass[12pt, reqno]{amsart}
\usepackage{amssymb}
\newtheorem{theorem}{Theorem}
\newtheorem*{theorem*}{Theorem}

\theoremstyle{remark}
\newtheorem*{definition*}{Definition}
\newtheorem*{remark*}{Remark}
\newtheorem{example}[theorem]{Example}
\newtheorem*{example*}{Example}

\begin{document}
\title[Non-Archimedean SMT for Hypersurfaces]{A Non-Archimedean Second Main Theorem for Hypersurfaces in Subgeneral
Position}

\date{October 29, 2024}
\begin{abstract}We apply an idea of Levin to obtain a non-truncated 
second main theorem for non-Archimedean analytic maps approximating 
algebraic hypersurfaces in subgeneral position. In some cases, for
example when all the hypersurfaces are non-linear and all the 
intersections are transverse, this
improves an inequality of Quang, whose inequality is sharp for the 
case of hyperplanes in subgeneral position.
\end{abstract}
\author{Ta Thi Hoai An}
\address{Institute of Mathematics\\
Vietnam Academy of Science and Technology\\
18 Hoang Quoc Viet Road, Cau Giay District\\
10307 Hanoi, VIETNAM}
\email{tthan@math.ac.vn}
\author{William Cherry}
\address{Department of Mathematics\\University of North Texas\\ 
1155 Union Circle \#311430, Denton, TX  76203\\USA}
\email{wcherry@unt.edu}
\author{Nguyen Viet Phuong}
\address{Thai Nguyen University of Economics and Business Administration, Vietnam} \email{nvphuongt@gmail.com}
\subjclass[2020]{30G06}
\keywords{non-Archimedean analytic curve, subgeneral position}
\thanks{The first and third authors are supported   by Institute of Mathematics, Vietnam Academy of Science and Technology  under grant number NVCSCL2023--2024 and Vietnam's National Foundation for Science and Technology Development (NAFOSTED) under grant number  101.04-2021.41.}
\maketitle
Ru \cite{Ru} observed that in non-Archimedean value distribution theory,
the second main theorem without truncation or ramification for
non-Archimedean analytic curves approximating algebraic 
hypersurfaces in projective
space follows from 
the first main theorem. An \cite{An} generalized Ru's approach to projective
hypersurfaces in general position with a projective variety. Recently,
Quang \cite{Quang} applied this approach to hypersurfaces
in subgeneral position. Although conjecturally the degree of the hypersurfaces
should come into these inequalities, in the works cited above, as with most
contemporary work in value distribution theory, only the number of
hypersurfaces and their intersection combinatorics enters into the 
inequalities. Prior work of Levin \cite{Levin} is one of the few cases where
the degrees of the hypersurfaces come into the inequality, although not in
what is conjecturally believed to be the optimal way. Levin considered only
hypersurfaces in general position. The purpose of this note is to record what
Levin's approach yields when the hypersurfaces are allowed to be in subgeneral
position. Levin's approach is most useful when the hypersurfaces meet
transversely. Although hypersurfaces in subgeneral position do not often 
meet transversely, when they do, the approach here sometimes 
gives an improvement on Quang's inequality.
\par
We'll take terminology and notation as in \cite[\S2]{Levin}.
\par
As in \cite{Levin}, we consider a projective variety $X\subseteq\mathbf{P}^N$
of dimension $n$
over a complete algebraically closed non-Archimedean field $K$ of arbitrary
characteristic. We consider $q$ projective hypersurfaces $D_1,\dots,D_q$
defined over $K$ in $\mathbf{P}^N.$ We do not necessarily assume that the
$D_j$ are distinct. We assume that $X$ is not entirely
contained in any of the hypersurfaces $D_j.$ For convenience,
we define the dimension of the empty set to be $-1.$ Then, for each integer
$m$ from $-1$ to $n-1,$ we define $t_m$ to be the smallest integer
such that for every subset $I\subseteq\{1,\dots,q\}$ of cardinality
$t_m+1,$ we have
$$
\dim\left(\bigcap_{i\in I}D_i\cap X\right)\le m.
$$
Note that if $\dim \left(X \cap D_1 \cap \dots \cap D_q\right)>m,$
we take $t_m=q.$ Most important for us will be the two numbers
$t_{-1}$ and $t_0.$ As each hypersurface can lower the dimension by at most
one, we have $t_m\ge\min\{n-m,q\}$ for each $m.$ When equality holds for all
$m,$ the hypersurfaces are said to be in general position with $X.$
More generally, the hypersurfaces are said to be in $t_{-1}$-subgeneral 
position with $X,$ and this has meaning when $q>t_{-1}.$
\par
Then, with this notation, Quang's subgeneral position without truncation result
is:

\begin{theorem}[{\cite[Th.\ 2]{Quang}}]\label{QuangThm}
Let $f:K\rightarrow X$ be a non-constant non-Archimedean analytic map
not entirely contained in any of the hypersurfaces $D_j.$
For $r>0,$
$$
\sum_{j=1}^q\frac{m_f(r,D_j)}{\deg D_j}\le t_{-1} T_f(r)+O(1).
$$
\end{theorem}
When $X=\mathbf{P}^N$ and the $D_j$ are all hyperplanes $H_j,$
then one can easily see that Quang's result is best possible in that
the coefficient $t_{-1}$ cannot be reduced.

\begin{example}\label{QuangSharp}Let $H_1,\dots,H_q$ be any collection of not necessarily 
distinct hyperplanes in $\mathbf{P}^N.$ Let $t_{-1}$ be defined as above,
namely the smallest integer such that the intersection of any
$t_{-1}+1$ of the hyperplanes is empty (or $q$ if no intersection of the
hyperplanes is empty). Assume that the absolute value on $K$ is non-trivial,
meaning that there is at least one element of $K$ with positive
absolute value less than one. Then, there is an algebraically
non-degenerate non-Archimedean analytic map $f:K\rightarrow\mathbf{P}^N$
such that for all $r>0,$
$$
\sum_{j=1}^q m_f(r,H_j)\ge t_{-1}T_f(r)-o(T_f(r)),
$$
and so the coefficient $t_{-1}$ in Quang's inequality
cannot be reduced.
\end{example}
\begin{proof}
By the definition of $t_{-1},$
we assume, without loss of generality, that $H_1,\dots,H_{t_{-1}}$ each pass
through the point in projective space with coordinates $(0,0,\dots,0,1).$
Let $f_1,\dots,f_N$ be any transcendental entire functions on $K$
such that $T_{f_i}(r)=o(T_{f_j}(r))$ for all $1\le i < j \le N.$
For example, since we assumed that the absolute value on $K$ is non-trivial,
let $a$ be an element of $K$ such that $0<|a|<1.$
Then, 
$$
g(z)=\sum_{n=0}^\infty(a^nz)^n
$$
has infinite radius of convergence and hence is a transcendental entire
function on $K.$ Now let $f_0=1,$ $f_1=g,$ $f_2=g\circ g,$
$f_3=g\circ g\circ g,$
and so on. Then, $T_{f_i}(r)=o(T_{f_j}(r))$ for all $1\le i < j \le N$
by \cite[Th.~2.44]{HuYang}. Now, let
$f=(f_0,f_1,,f_2,\dots,f_{N-1},f_N),$ and observe that
$T_f(r)=T_{f_N}(r)+o(T_{f_N}(r)).$
Let $P(X_0,\dots,X_N)$ be any non-trivial
homogeneous polynomial. Let $k$ be the largest integer such that $X_k$
appears in a monomial in $P$ with non-zero coefficient. If $k=0,$
then trivially, $P\circ f\not\equiv0.$ If $k>0,$ then
let $d$ be the $X_k$-degree of $P.$ Then, 
$T_{P\circ f}(r)=dT_{f_k}(r)+o(T_{f_k}(r))$ by \cite[Th.~2.11]{HuYang},
and so $P\circ f\not\equiv 0,$ and thus $f$ is algebraically non-degenerate.
\par
Because $H_1,\dots,H_{t_{-1}}$ contain the point $(0,0,\dots,0,1),$
we have that 
$$m_f(r,H_j)\ge T_{f_N}(r)-o(T_{f_N}(r))
\qquad\textrm{for all}~j~\textrm{from}~1~\textrm{to}~t_{-1}.\qedhere$$ 
\end{proof}
\par
\begin{remark*}If instead, the absolute value on $K$ is trivial, then
the only entire functions on $K$ are the polynomials, and all
analytic maps from $K$ to $\mathbf{P}^N$ are algebraic. Still, with
hyperplanes as in Example~\ref{QuangSharp} and given $\varepsilon>0,$
one easily sees that for all large enough integers $d,$ depending on
$\varepsilon,$
$$
f(z)=(1,z,z^2,\dots,z^{N-1},z^d)
$$
is a linearly non-degenerate map $f:K\rightarrow \mathbf{P}^N$ such
that for all $r$ sufficiently large,
$$
\sum_{j=1}^q m_f(r,H_j)\ge (t_{-1}-\varepsilon)T_f(r).
$$
\end{remark*}
\par
One expects to be able to do better when the $D_j$ are non-linear
hypersurfaces, but to date, there are few ideas about how to take
advantage of higher degree $D_j$ to improve this type of inequality.
Levin showed one way to sometimes be able to take advantage of higher degree,
particularly when the $D_j$ intersect transversely:

\begin{theorem}[{\cite[Th.\ 10]{Levin}}]\label{LevinThm}
Let $X\subseteq\mathbf{P}^N$ be a projective variety over $K$ of dimension
$n\ge1.$ Let $D_1,\dots,D_q$ be hypersurfaces in $\mathbf{P}^N$ over $K$
that are in general position with $X.$ Let $M$ be the smallest positive
integer such that for any subset $I\subseteq\{1,\dots,q\}$ of cardinality $n,$
$$
\bigcap_{i\in I}D_i\cap X \subseteq M \mathrm{Supp}\left(\bigcap_{i \in I}
D_i\cap X\right),
$$
where we view $D_1,\dots,D_q$ and $X$ as closed subschemes of $\mathbf{P}^N.$
Let \hbox{$f:K\rightarrow X$}
be a nonconstant non-Archimedean analytic map whose
image is not completely contained in any of the hypersurfaces $D_j.$
Then, for all $r\ge1,$
$$
\sum_{j=1}^q\frac{m_f(r,D_j)}{\deg D_j}\le\left(n-1+\max_{1\le j \le q}
\frac{M}{\deg D_j}\right)T_f(r)+O(1).
$$
Note that if all the intersections amongst $X$ and $D_1,\dots,D_q$ are
transverse, then $M=1.$
\end{theorem}
If we relax the general position hypothesis, then Levin's idea gives
\begin{theorem}\label{NewThm}
Let $X\subseteq\mathbf{P}^N$ be a projective variety over $K$ of dimension
$n\ge1.$ Let $D_1,\dots,D_q$ be not necessarily distinct
hypersurfaces in $\mathbf{P}^N$ over $K.$ Let $t_0$ and $t_{-1}$
be defined as above.
Let $M$ be the smallest positive
integer such that for any subset $I\subseteq\{1,\dots,q\}$ of cardinality
$t_0+1,$
$$
\bigcap_{i\in I}D_i\cap X \subseteq M \mathrm{Supp}\left(\bigcap_{i \in I}
D_i\cap X\right),
$$
where we view $D_1,\dots,D_q$ and $X$ as closed subschemes of $\mathbf{P}^N.$
Let \hbox{$f:K\rightarrow X$} be a nonconstant non-Archimedean analytic map whose
image is not completely contained in any of the hypersurfaces $D_j.$
Let
$$
\alpha=\max_{\begin{array}{c} I\subseteq\{1,\dots,q\}\\ |I|=t_{-1}-t_0
\end{array}}
\sum_{i\in I}\min\left\{\frac{M}{\deg D_i},1\right\}.
$$
Then, for all $r\ge1,$
$$
\sum_{j=1}^q\frac{m_f(r,D_j)}{\deg D_j}\le (t_0+\alpha)T_f(r)+O(1).
$$
\end{theorem}
Note that if all the intersections amongst $X$ and $D_1,\dots,D_q$ are
transverse, then $M=1.$
\begin{remark*}In case that the $D_j$ are in general position and $q>n,$ then
\hbox{$t_0=n-1,$} $t_{-1}=n,$ $t_{-1}-t_0=1,$ and
$$
\alpha=\max_{1\le j \le q}\min\left\{\frac{M}{\deg D_j},1\right\},
$$
and so we recover Theorem~\ref{LevinThm}. Moreover, $t_0+\alpha\le t_{-1},$
and so we also recover Theorem~\ref{QuangThm}.
\par
\begin{remark*}Here we, as was Levin in \cite{Levin}, are mainly only
interested in the case that none of the hypersurfaces $D_i$ are hyperplanes.
In Levin's formulation, if any of the $D_i$ are hyperplanes, then
Levin's inequality does not improve upon \cite{An}.
In a recent extension of
Levin's work in a different direction,
Huynh \cite{Tuan} shows that Levin's inequality can 
sometimes be improved in the case where $X=\mathbf{P}^N,$ when $q=N,$
and when the degeneracy of $f$ is further restricted by not allowing
the image of $f$ to be contained in certain hyperplanes tangent to the
hypersurfaces $D_i,$
which then can give a better result when some of the hypersurfaces are 
hyperplanes, for example in the case of a non-singular conic and a line
meeting transversally in $\mathbf{P}^2.$ Combining Huynh's idea
with ours might result in some improvement to our inequality under
the additional assumption that the image of $f$ is not contained in certain
tangent hyperplanes to the hypersurfaces $D_j,$ but we do not now see a
nice formulation of such a result in a case more general than what Huynh
already treats.
\end{remark*}
\par
Although not typical for hypersurface arrangements not in
general position, it is possible for hypersurfaces not in general position
to intersect transversely when all the extra intersections happen in
dimension zero, for example three conics in $\mathbf{P}^2$ all intersecting
transversely in a common set of points. This is precisely the case where
Levin's idea gives the best improvement.  In the case of three conics
intersecting in $\mathbf{P}^2$ transversely in a common set of points,
we have $t_0=1$ and $t_{-1}=q=3,$ and $\alpha=1.$ In this case,
Theorem~\ref{QuangThm}
gives
$$
\sum_{j=1}^3\frac{m_f(r,D_j)}{2}\le 3T_f(r)+O(1),
$$
which is just the First Main Theorem, but Theorem~\ref{NewThm} gives the
better
$$
\sum_{j=1}^3\frac{m_f(r,D_j)}{2}\le2T_f(r)+O(1).
$$
\end{remark*}
\par
For some collections of non-linear hypersurfaces,
Theorem~\ref{NewThm} is a sharp improvement of Theorem~\ref{QuangThm}.
\begin{example}
Let $K$ be any algebraically closed complete non-Archimedean field with
characteristic different from three. Let 
\begin{eqnarray*}
Q_1(X_0,X_1,X_2)&=&X_0X_1-X_2^2,\\
Q_2(X_0,X_1,X_2)&=&X_0X_2-X_1^2,\\
\textrm{and}~Q_3(X_0,X_1,X_2)&=&Q_1(X_0,X_1,X_2)+3Q_2(X_0,X_1,X_2).
\end{eqnarray*}
For $j=1,2,$ and $3,$ let $D_j$ be the hypersurface in $\mathbf{P}^2$
determined by the vanishing of $Q_j.$ Then, the $D_j$ are three conics
in $\mathbf{P}^2$ intersecting transversely in the four points
$(1,0,0),$ $(1,1,1),$ $(1,\omega,\omega^2),$ and $(1,\omega^2,\omega),$
where $\omega$ is any primitive third root of unity.
Let $f:K\rightarrow\mathbf{P}^2$ be the non-constant linear map
with coordinate functions $(z,1,0).$ Then, for $r\ge1,$
\begin{eqnarray*}
m_f(r,D_1)&=&\log\frac{|z|_r^2}{|z|_r}=\log r\\
m_f(r,D_2)&=&\log\frac{|z|_r^2}{|-1|_r}=2\log r\\
\textrm{and}~m_f(r,D_3)&=&\log\frac{|z|_r^2}{|z-3|_r}=\log r.
\end{eqnarray*}
Hence, for $r\ge1,$
$$
\sum_{j=1}^3 \frac{m_f(r,D_j)}{2}=2\log r = 2T_f(r),
$$
showing that Theorem~\ref{NewThm} is sharp in this case.
\end{example}
\begin{proof}[Proof of Theorem~\ref{NewThm}]
We essentially follow \cite[Th.\ 10]{Levin} but allow that taking 
additional hypersurfaces when the intersection has dimension zero may
not further reduce the dimension.
\par
Let $f=(f_0,\dots,f_N),$ where $f_0,\dots,f_N$ are entire
without common zeros. Let $D_1,\dots,D_q$ be defined by
homogeneous polynomials
$Q_1,\dots,Q_q$ in $K[X_0,\dots,X_N].$
\par
If there are fewer than $t_0$ indices $j\in\{1,\dots,q\}$ such that
$m_f(r,D_j)\ne O(1)$ as $r\to\infty,$ then the theorem follows from 
the First Main Theorem. We therefore assume from now on that there
are at least $t_0$ indices $j$ such that $m_f(r,D_j)\to\infty$
as $r\to\infty.$
\par
Let $I,J\subseteq\{1,\dots,q\}$ be such that $m_f(r,D_i)=O(1)$
for all $i\in I,$ and such that $m_f(r,D_j)\to\infty$ as $r\to\infty$
for all $j\in J.$ Let $r_0$ be large enough so that
for all $i\in I,$ all $j\in J,$ and all $r\ge r_0,$
$$
\frac{m_f(r,D_j)}{\deg D_j} \ge \frac{m_f(r,D_i)}{\deg D_i}
$$
Fix $r\ge\max\{1,r_0\}.$ After reindexing, we may assume that
$$
\frac{m_f(r,D_1)}{\deg D_1}\ge \frac{m_f(r,D_2)}{\deg D_2}
\ge\dots\ge \frac{m_f(r,D_q)}{\deg D_q}.
$$
If $t_{-1}<q,$ then 
$$
D_1\cap\dots\cap D_{t_{-1}}\cap D_{t_{-1}+1}=\emptyset.
$$
Hilbert's Nullstellensatz then tells us that for each of the coordinate
functions $X_j,$ there is some power $m_j$ and some homogeneous
polynomials $A_{ji}$ such that
$$
X_j^{m_j}=\sum_{i=1}^{t_{-1}+1}A_{ji}Q_i.
$$
Of course, $\deg A_{ji}=m_j-\deg D_i.$ Thus, there exists a constant
$C_1,$ depending only on the polynomials $A_{ji},$ such that for all
$j=0,\dots,N,$
$$
|f_j|_r^{m_j}\le C_1\max_{1\le j \le t_{-1}+1}|f|_r^{m_j-\deg D_i}
|Q_i\circ f|_r.
$$
Choosing $j$ so that $|f_j|_r=|f|_r$ and canceling $|f|_r^{m_j}$
from both sides then gives
$$
1\le C_1\max_{1\le i \le t_{-1}+1}\frac{|Q_i\circ f|_r}{|f|_r^{\deg D_i}}\le
C_1\frac{|Q_{t_{-1}+1}\circ f|_r}{|f|_r^{\deg D_{t_{-1}+1}}}.
$$
Hence,
\begin{equation}\label{lastpiece}
\frac{m_f(r,D_i)}{\deg D_i}\le\frac{\log C_1}{\deg D_i},\textrm{~for all~} i>t_{-1}.
\end{equation}
Note that the constant
$C_1$ was chosen depending on the $A_{ji},$ which depend on $r$ in the
sense that they depend on how we indexed $Q_1,\dots,Q_q.$ As there are
only finitely many ways of reindexing, $C_1$ can be chosen independent of $r.$
\par
If $t_0=t_{-1}=q,$ then the theorem follows from the First Main Theorem,
or from Theorem~\ref{QuangThm}, so we henceforth assume that $t_0<q.$
Levin's idea allows us to handle $m_f(r,D_j)$ for
$j=t_0+1,\dots,t_{-1}.$ Fix one such index $j.$ By the definition of $t_0,$
$$
X\cap D_1 \cap \dots \cap D_{t_0} \cap D_j
$$
is a finite set of points $\{P_1,\dots,P_s\}.$ We now want to
choose hyperplanes
$H_1,\dots,H_s$ so that the following conditions are satisfied:
\begin{itemize}
\item $P_i \in H_i$ for $i=1,\dots,s;$
\item $P_j \notin H_i$ for all $ j\ne i\in \{1,\dots,s\};$
\item If we define $E_i=D_i$ for $i=1,\dots,t_0$ and
$E_i=H_{i-t_0}$ for $i=t_0+1,\dots,t_0+s,$ then for any index set
$I\subseteq\{1,\dots,t_0+s\}$ with cardinality $|I|\ge t_0+2,$
$$
X\cap\bigcap_{i\in I}E_i=\emptyset;
$$
\item The image of $f$ is not completely contained in any of the $H_i.$
\end{itemize}

We now re-order the $E_i$ so that
$$
\frac{m_f(r,E_1)}{\deg E_1}\ge\dots\ge \frac{m_f(r,E_{t_0+s})}{\deg E_{t_0+s}}.
$$
By our choice of the hyperplanes, if $s\ge2,$ then 
$$
X\cap E_1 \cap \dots \cap E_{t_0} \cap E_{t_0+1} \cap E_{t_0+2}=\emptyset.
$$
As before, we apply the Nullstellensatz, to find a constant $C_2$ such 
that for all $i>t_0+1,$
\begin{equation}\label{finalstep}
\frac{m_f(r,E_i)}{\deg E_i}\le \frac{\log C_2}{\deg E_i},
\end{equation}
where $C_2$ appears to depend on $r$ and our choice of the index $j,$
as it depends on
$\{P_1,\dots,P_s\}$ and the choice of hyperplanes $H_1,\dots,H_s,$
and these depend on the index $j$ and on our original re-indexing
of $D_1,\dots,D_q.$ 
But again, there are only finitely many ways to re-index $D_1,\dots,D_q,$
and for each reindexing, there are only finitely many choices of the index
$j,$ and so $C_2$ can be taken to be independent of $r.$
\par
By enlarging $r_0$ if necessary and by our assumption that there were
at least $t_0$ indices $j$ such that $m_f(r,D_j)\to\infty$ as $r\to\infty,$
we conclude that all the $E_i$ in~(\ref{finalstep}) are, in fact, hyperplanes.
Hence,
\begin{equation}\label{hypbound}
\sum_{i=1}^s m_f(H_i,r)\le T_f(r)+O(1).
\end{equation}
\par
Now, choose linear defining forms $L_i$ for each of the hyperplanes $H_i.$
By our definition of $M,$
$$
X\cap D_1 \cap \dots \cap D_{t_0} \cap D_j\subseteq
M\mathrm{Supp}\left(
X\cap D_1 \cap \dots \cap D_{t_0} \cap D_j\right).
$$
Thus, there exist homogeneous polynomials $B_1,\dots,B_{t_0}$ and $B_j$
with $$\hbox{$\deg B_i=Ms-\deg D_i$} $$ such that
$$
(L_1 \cdots L_s)^M-B_jQ_j-\sum_{i=1}^{t_0}B_iQ_i
$$
vanishes on $X.$ Therefore,  there is a constant $C_3,$ depending only on 
the polynomials $B_i,$ and hence can be taken to be independent of $r$
if we consider all possible reindexings of $Q_1,\dots,Q_q,$ such that
$$
\prod_{i=1}^s |L_i\circ f|_r^M \le C_3
\max_{k\in\{1,\dots,t_0\}\cup\{j\}}|f|_r^{Ms-\deg D_k}
|Q_k\circ f|_r \le C_3|f|_r^{Ms}
\frac{|Q_j\circ f|_r}{|f|_r^{\deg D_j}}.
$$
Dividing both sides by $|f|_r^{Ms}$ and taking logarithms, we get
$$
m_f(r,D_j) \le
M\sum_{i=1}^s m_f(r,H_i)+C_4,
$$
with $C_4$ independent of $r,$ again by observing that there are only
finitely many ways to reindex the $D_j.$ Combining with $(\ref{hypbound}),$
we conclude
$$
\frac{m_f(r,D_j)}{\deg D_j} \le \frac{M}{\deg D_j}T_f(r)+O(1).
$$
Note that, in any case,
$$
\frac{m_f(r,D_j)}{\deg D_j}\le T_f(r)+O(1)
$$
by the First Main Theorem, and so we in fact have
\begin{equation}\label{middlepiece}
\frac{m_f(r,D_j)}{\deg D_j}\le \min
\left\{\frac{M}{\deg D_j},1\right\}T_f(r)+O(1).
\end{equation}
\par
We now split the left-hand-side of the inequality we want to prove into
three pieces:
$$
\sum_{j=1}^q\frac{m_f(r,D_j)}{\deg D_j} \le
\sum_{j=1}^{t_0}\frac{m_f(r,D_j)}{\deg D_j}
+\sum_{j=t_0+1}^{t_{-1}}\frac{m_f(r,D_j)}{\deg D_j} +
\sum_{j=t_{-1}+1}^q\frac{m_f(r,D_j)}{\deg D_j}.
$$
Using~(\ref{lastpiece}), we can replace the last piece with $O(1).$
Now, by our choice of ordering for the $D_j$
and~(\ref{middlepiece}), we can reduce the middle piece
to
$$
\sum_{j=t_0+1}^{t_{-1}}\frac{m_f(r,D_j)}{\deg D_j}\le
\alpha T_f(r)+O(1),
$$
where
$$
\alpha=\max_{\begin{array}{c} I\subseteq\{1,\dots,q\}\\ |I|=t_{-1}-t_0
\end{array}}
\sum_{i\in I}\min\left\{\frac{M}{\deg D_i},1\right\}.
$$
Putting those together and applying the First Main Theorem one final time,
we now have
\begin{align*}
\sum_{j=1}^q\frac{m_f(r,D_j)}{\deg D_j}&\le
\sum_{j=1}^{t_0}\frac{m_f(r,D_j)}{\deg D_j}+
\alpha T_f(r)+O(1)\\
&\le t_0T_f(r)+\alpha T_f(r)+O(1)\\
&\le (t_0+\alpha )T_f(r)+O(1). \qedhere
\end{align*}
\end{proof}

\section*{Acknowledgment}
 A part of this article was written while the first and second  authors were
visiting the Vietnam Institute for Advanced Study in Mathematics (VIASM). They
would like to thank the institute for warm hospitality and  support. 
The authors greatly thank the referee for pointing out a problem with 
Example~\ref{QuangSharp} in an earlier draft and for suggesting
inequality~(\ref{hypbound}) and the addition of the minimum
in~(\ref{middlepiece}). These suggestions both significantly simplified and 
improved our original formulation of Theorem~\ref{NewThm}.


\begin{thebibliography}{[MMM]}
\bibitem[A]{An} T.~T.~H.~An, A Defect Relation for Non-Archimedean Analytic
Curves in Arbitrary Projective Varieties, Proc.\ Amer.\ Math.\ Soc.\
135 (2007), 1255--1261.
\bibitem[HY]{HuYang} P.-C.~Hu and C.-C.~Yang, \textit{Meromorphic
Functions over Non-Archimedean Fields,} Kluwer, 2000.
\bibitem[H]{Tuan} D.~T.~Huynh, A $p$-adic Second Main Theorem,
Preprint, 2024, arXiv: 2405.14197
\bibitem[L]{Levin} A.~Levin, On the $p$-Adic Second Main Theorem,
Proc.\ Amer.\ Math.\ Soc.\ 143 (2015), 633--640.
\bibitem[Q]{Quang} S.~D.~Quang, Truncated Second Main Theorem for
Non-Archimedean Meromorphic Maps, Preprint, 2023, arXiv: 2306.07594
\bibitem[R]{Ru} M.~Ru, A Note on $p$-adic Nevanlinna Theory, 
Proc.\ Amer.\ Math.\ Soc.\ 6 (1995), 1263--1269.
\end{thebibliography}
\end{document}